\documentclass[]{siamart1116}

\usepackage{lipsum}
\usepackage{amsfonts}
\usepackage{graphicx}
\usepackage{epstopdf}
\usepackage{algorithmic}

\usepackage{amsmath}
\usepackage{amssymb,dsfont} 
\usepackage{graphicx}
\usepackage{booktabs}
\ifpdf
  \DeclareGraphicsExtensions{.eps,.pdf,.png,.jpg}
\else
  \DeclareGraphicsExtensions{.eps}
\fi

\numberwithin{theorem}{section}

\newcommand{\TheTitle}{How to avoid the curse of dimensionality: scalability of particle filters with and without importance weights} 
\newcommand{\TheAuthors}{S.C. Surace, A. Kutschireiter, and J.-P. Pfister}

\headers{Curse of dimensionality in particle filters}{\TheAuthors}

\title{{\TheTitle}\thanks{Submitted to the editors Tuesday, April 11, 2017.
\funding{This work was funded by the Swiss National Science Foundation, grants PP00P3 150637 and PZ00P3 137200.}}}

\author{
 Simone Carlo Surace\thanks{Institute of Neuroinformatics, University of Z\"urich and ETH Z\"urich, Z\"urich, Switzerland.
    (\email{surace@ini.uzh.ch)}}
  \and
  Anna Kutschireiter\footnotemark[2]
  \and
  Jean-Pascal Pfister\footnotemark[2]
}

\usepackage{amsopn}

\newcommand{\parenths}[1]{\left( #1 \right)}
\newcommand{\brackets}[1]{\left[ #1 \right]}
\newcommand{\braces}[1]{\left\{ #1 \right\}}

\newcommand{\norm}[1]{\left\Vert #1 \right\Vert}

\newcommand{\com}[1]{{\vspace{5mm} \noindent\Large\color{red}#1\color{black}\vspace{5mm}}}

\newcommand{\soutt}[1]{\sout{#1} }
\newcommand{\temporary}[1]{#1}
\newcommand{\EE}[1]{\mathds{E}\brackets{#1}}
\newcommand{\dbar}[1]{\overline{\overline{#1}}}

\renewcommand{\soutt}[1]{}   
\renewcommand{\com}[1]{} 
\renewcommand{\temporary}[1]{} 

\ifpdf
\hypersetup{
  pdftitle={\TheTitle},
  pdfauthor={\TheAuthors}
}
\fi

\begin{document}

\maketitle

\begin{abstract}
Particle filters are a popular and flexible class of numerical algorithms to solve a large class of nonlinear filtering problems. 
However, standard particle filters with importance weights have been shown to require a sample size that increases exponentially with the dimension $D$ of the state space in order to achieve a certain performance, which precludes their use in very high-dimensional filtering problems. 
Here, we focus on the dynamic aspect of this Ôcurse of dimensionalityÕ (COD) in continuous time filtering, which is caused by the degeneracy of importance weights over time. 
We show that the degeneracy occurs on a time-scale that decreases with increasing $D$. 
In order to soften the effects of weight degeneracy, most particle filters use particle resampling and improved proposal functions for the particle motion. 
We explain why neither of the two can prevent the COD in general. 
In order to address this fundamental problem, we investigate an existing filtering algorithm based on optimal feedback control that sidesteps the use of importance weights. 
We use numerical experiments to show that this Feedback Particle Filter (FPF) by \cite{Yang:hf} does not exhibit a COD. 
\end{abstract}



\section{Introduction}
Filtering problems are rarely exactly solvable with a finite amount of computational resources, requiring numerical techniques in order to approximately represent or sample from the filtering distribution. 
Particle filters, which have been first introduced by \cite{Gordon:1993be}, have seen widespread use as general-purpose algorithms to solve nonlinear filtering problems (see \cite{Doucet:2009us,Kunsch:2013gg,Crisan:2002hq} for a general survey). 
They use sequential importance sampling in order to calculate the filtering distribution. 
At each time-step, samples or particles are drawn from a proposal density (such as the prior transition probability) and then re-weighted according to the observations. 
Particle filters have very few restrictions with regards to the type of generative model that underlies the filtering problem;
they can be applied to highly nonlinear and non-Gaussian models, which is an advantage compared to the well-known Kalman filters.

However, they suffer from a very severe problem that affects all importance sampling-based algorithms; their efficiency diminishes rapidly with increasing dimension of the state space. 
This `curse of dimensionality' (COD) occurs because in high-dimensional spaces the importance weights are more likely to be degenerate, i.e. only a few weights are significant and all others are very close to zero.
There have been numerous studies on the COD in particle filtering \cite{Daum:2003tf,Bengtsson:2008ei,Bickel:2008bx,Snyder:2008by}, which will be described below in more detail.

In this paper, we first revisit the COD in particle filters with a focus on the dynamic aspect of the problem.
We study how the time-scales of weight degeneracy and of the performance benefit due to resampling scale with dimension.
We find that the time-scale of weight degeneracy is inversely proportional to the dimensionality and proportional to the logarithm of ensemble size.
From this, we obtain an exponential scaling of the ensemble size that is required to obtain a fixed time-scale of weight degeneracy.
We argue that because the resampling-induced benefits occur at a time-scale that scales weakly with dimension, this exponential scaling is necessary in order to avoid the collapse of the particle filter.

Next, we show that the COD can be avoided by particle filters that do not employ importance weights.
Such a filter must sample the posterior distribution directly, and the particle motions have to fully take into account the observations.
An unweighted particle filter for the classical filtering problem, the Feedback Particle Filter (FPF) by \cite{Yang:hf}, uses optimal feedback control in order to guide the particle motions according to the observations.
As a result, the particles of the FPF are samples from the posterior distribution.
Although the authors claim that the FPF method does not suffer from the COD, to our knowledge this has never been explicitly demonstrated.
We fill this gap in the literature by demonstrating that the FPF does not suffer from the COD.

\subsection{Existing literature on the COD in particle filtering}
The COD of particle filters has been studied using simple back-of-the-envelope calculations \cite{Daum:2003tf,Snyder:2008by} and by proving statements on the convergence of the maximum weight \cite{Bengtsson:2008ei,Bickel:2008bx}.
These works reduce the filtering problem to a single Bayesian update step.
In doing so, the COD of particle filters is reduced to the general COD of importance sampling and the effects of proposal distributions or resampling on the COD are not explored in detail.

The effect of proposal distributions on particle filtering in high-dimensional problems has received some attention.
In discrete time, particle filters with `smart' proposal distributions (i.e. with improved laws that govern the motion of particles) have been shown to perform well in certain high-dimensional problems \cite{vanLeeuwen:2009im,vanLeeuwen:2010ho,Ades:2014gn}, leading some authors to conjecture that the COD could be avoided using carefully crafted proposal distributions. 
However, it has been argued \cite{Snyder:2015dw} that the system considered in \cite{Ades:2014gn} (and similar systems) was effectively a low-dimensional system disguised as a high-dimensional one.
A precise definition of a sequence of filtering problems of increasing dimension therefore requires a notion of `effective dimension', which also plays a role in the arguments brought forward by \cite{Daum:2003tf} and \cite{Snyder:2008by}.
Special properties of the model sequence, e.g. a low-dimensional dynamical system that is embedded in spaces of increasing dimensions, could potentially avoid the COD as a function of $D$, but not as a function of effective dimension.
Improved proposals are usually designed to reduce weight degeneracy by minimizing the rate of change of the variance of the importance weights.
But even the optimal proposal function is not able to completely prevent the weight degeneracy.
Moreover, in continuous time it is not known whether an optimal proposal even exists.

Meanwhile, the existing literature offers only a cursory treatment of the effect of resampling on the dimensionality-dependent scaling of particle filters. 
While it is understood that the most widespread resampling algorithm, multinomial resampling, is suboptimal compared to branching (see \cite{Crisan:1998un}, \cite{Chopin:2004cn}, \cite{Douc:2005wa}, \cite{Crisan:2006dn}, and also \cite{Bain:2009id}, p.250), we cannot be certain about the existence of a resampling scheme that resolves the COD.
The arguments against such a possibility, which are brought forward e.g. in \cite{Snyder:2015dw} and \cite{Chopin:2004cn}, are mostly heuristic and can be summarized as follows:
while resampling resets the particle weights to undo the effects of weight degeneracy, it does not improve the quality of the sample instantaneously.
The potential improvement is temporary and comes from particle motions that originate from regions of high likelihood.
In the next section we will look more closely at this mechanism and thereby support the explanation provided by \cite{Snyder:2015dw}.

The conjecture that the COD may be avoided by particle filters without importance weights has been brought up in previous works.
For example, the (stochastic) particle flow filter \cite{Daum:2010bu,DeMelo:2015vz} holds the promise of avoiding the COD \cite{Daum:2017ua}.
This filter implements an implicit Bayesian update by defining a flow field at each measurement:
whenever a measurement is registered, a little loop in synthetic time is inserted, during which the particles evolve according to this flow field, and are thus propagated from the prediction to the posterior density.
A similar approach, the continuous-discrete particle filter \cite{Yang:2014cv}, also relies on inserting an artificial loop at the times of measurement arrival, but the particles are propagated according to an innovation-error based gain feedback structure.
Since particle flow filters require a clear separation between prediction and update in their respective filtering algorithms, it is not clear how the results would transfer to continuous-time measurement models, where prediction and update have to be processed simultaneously.

\section{Preliminaries}
\label{problem}
We restrict our overall discussion to the classical nonlinear filtering problem in continuous time, where the state $X_t\in\mathds{R}^D$ and the observation $Y_t\in\mathds{R}^D$ are $D$-dimensional diffusion processes that are solutions to the It\^o stochastic differential equations (SDE)
\begin{align}
dX_t &= f(X_t)dt+g(X_t)dW_t, \\
dY_t &= h(X_t)dt+dV_t,
\end{align}
where $W_t$ and $V_t$ are independent $D$-dimensional Brownian motions.
The stochastic filtering problem is to find conditional expectations $\mathds{E}\brackets{\varphi(X_t)|\mathcal{F}_t^Y}$, where $\varphi:\mathds{R}^D~\rightarrow~\mathds{R}$ is a measurable function and $\mathcal{F}_t^Y$ is the filtration generated by the observation process. 
Next, we will describe the two approximate filtering algorithms that will be studied and compared in this paper.

\subsection{The Bootstrap Particle Filter (BPF)}
The bootstrap or vanilla particle filter uses importance sampling (IS) to approximate the conditional expectation as
\begin{equation}
\mathds{E}\brackets{\varphi(X_t)|\mathcal{F}_t^Y}     \approx      \dbar\varphi_t\doteq\sum_{i=1}^Nm_t^{(i)}\varphi(Z_t^{(i)}),
\label{Eqn:IS}
\end{equation}
where $Z_t^{(i)}\in\mathds{R}^D$ are the samples or `particles' that evolve according to the dynamics of the hidden state, and $m_t^{(i)}$ are the (normalized) importance weights. 
This particle filter is usually formulated in discrete time \cite{Gordon:1993be,Doucet:2009us}, but we will use the standard formulation of continuous-time particle filtering (see Ch.~9 of \cite{Bain:2009id} and Ch.~23 of \cite{Crisan:2011tua}).
Between resampling times, the time evolution of the particle system is given by
\begin{align}
\label{Eqn:particlemotion} dZ_t^{(i)} &= f(Z_t^{(i)})dt+g(Z_t^{(i)})dB_t^{(i)}, \\
\label{Eqn:weightupdate} m_t^{(i)}&=\frac{M_t^{(i)}}{\sum_{j=1}^NM_t^{(j)}}, \quad dM_t^{(i)}=M_t^{(i)}h(Z_t^{(i)})\cdot dY_t,
\end{align}
where $B_t^{(i)}$ are independent Brownian motions and $\cdot$ denotes the standard scalar product on $\mathds{R}^D$. 
In this formulation, the unnormalized importance weight $M_t^{(i)}$ arises as a Radon-Nikodym derivative $d\mathds{P}^{(i)}/d\tilde{\mathds{P}}$ of the measure $\mathds{P}^{(i)}$ under which the observations are generated from the state $Z_t^{(i)}$ to a reference measure $\tilde{\mathds{P}}$ under which $Y_t$ is a Brownian motion.

It is well established \cite{Doucet:2009us} that the BPF suffers from weight degeneracy, i.e. the weights evolve to become less equal, with only a few significant weights and all other weights being negligibly small. 
As a result, the variance of the IS estimate in Eq.~\eqref{Eqn:IS} grows and performance drops. 
It is useful to introduce an `effective sample size' $N_{\text{eff}}$ that measures the number of samples that would be required to match the performance of IS with a MC sampler that samples directly from the filtering distribution. 
A useful approximation of effective sample size is given by the inverse of the sum of squared weights, 
\begin{equation}
N_{\text{eff},t}\approx \brackets{\sum_{i=1}^N\parenths{m_t^{(i)}}^2}^{-1}\doteq \tilde N_{\text{eff},t},
\label{Eqn:ESS}
\end{equation}
(see \cite{Martino:2017bv} and references therein). 
As a measure that only depends on the importance weights and not on the locations of the samples, it is not a good predictor of performance when samples are clustered. 
In particular, the typical outcome of multinomial resampling from an impoverished sample is that there are only a few distinct particle positions and multiple particles at each of those positions. 
The new sample is still impoverished and has a large variance despite the fact that the above definition of $N_{\text{eff},t}$ yields a value of $N$ (after resampling, all weights are reset to a value of $1/N$). 
With this caveat in mind, we assume that $\tilde N_{\text{eff},t}$ in Eq.~\eqref{Eqn:ESS} gives an upper bound on the true value of $N_{\text{eff},t}$. 

Because of the weight degeneracy outlined above, the particle system has to be frequently resampled.
Resampling is performed according to a schedule that is usually based either on regular resampling intervals or on the effective sample size.
We use the convention that particle trajectories are right-continuous with left-side limits $Z_{t_r^-}^{(i)}$ encoding the position before resampling.
At each resampling time $t=t_r$, new particle positions $(\tilde Z_{t_r}^{(i)})_{i=1}^n$ are drawn from the set of positions $\{Z_{t_r^-}^{(i)},\, i=1,...,n\}$ with replacement.
The probability that $Z_{t_r^-}^{(i)}$ appears $n_i$ times within the tuple $(\tilde Z_{t_r}^{(i)})_{i=1}^n$ is chosen according to the multinomial distribution,
\begin{equation}
\text{Prob}\parenths{n_1,n_2,...,n_N}=N!\prod_{i=1}^N\frac{\parenths{m_t^{(i)}}^{n_i}}{n_i!}, \quad  \sum_{i=1}^Nn_i=N.
\label{Eqn:multinomial}
\end{equation}
The new position at $t=t_r$ is then set to the resampled position, i.e. $Z_{t_r}^{(i)}\leftarrow\tilde Z_{t_r}^{(i)}$ and all weights are reset to $1/N\leftarrow\tilde m_{t_r}^{(i)}$.
A discussion of various other forms of resampling can be found e.g. in \cite{Crisan:1998un}, \cite{Chopin:2004cn}, \cite{Douc:2005wa}, and \cite{Crisan:2006dn}.

\subsection{A particle filter without importance weights: the Feedback Particle Filter (FPF)}
In \cite{Yang:hf}, a particle filter was proposed that does not require importance weights, i.e. posterior expectations are approximated by unweighted averages,
\begin{equation}
\mathds{E}\brackets{\varphi(X_t)|\mathcal{F}_t^Y}     \approx      \bar{\varphi}_t\doteq\frac{1}{N}\sum_{i=1}^N\varphi(Z_t^{(i)})
\label{Eqn:MCS}
\end{equation}
(note the distinction to the weighted average $\dbar\varphi_t$ in Eq.~\eqref{Eqn:IS}).
By definition, the time evolution of particles fully incorporates the information given by the history of observations.
The particle system of the FPF evolves according to the It\^o SDEs
\begin{equation}
\label{Eqn:FPFparticles} \begin{split}
dZ_t^{(i)} &= \parenths{f(Z_t^{(i)})+\Omega(Z_t^{(i)},t)}dt\\
&\quad+g(Z_t^{(i)})dB_t^{(i)}+K(Z_t^{(i)},t)\brackets{dY_t-\tfrac{1}{2}(h(Z_t^{(i)})+\bar{h}_t)dt}, 
\end{split}
\end{equation}
where $K$ is a ($D\times D$ matrix-valued) function, and $\Omega$ is a $D$-dimensional vector-valued function with components given by
\begin{equation}
\Omega_i(z,t)=\frac{1}{2}\sum_{j=1}^D\sum_{k=1}^D K_{jk}(z,t)\frac{\partial}{\partial z_j}K_{ik}(z,t), \quad i=1,...,D.
\end{equation}
The function $K$ is the solution to an Euler-Lagrange boundary value problem (EL-BVP).
The latter results from an optimal control problem: 
the function $K$ is chosen such as to make the distribution of particle positions as close as possible to the posterior filtering distribution, and it ensures that the particle distribution converges to the true filtering distribution when $K$ is chosen according to the filtering distribution.
The function $K$, also called \emph{gain function}, is analogous to the Kalman gain of the Kalman-Bucy filter for the linear filtering problem.
As such, $K$ (and therefore $\Omega$) introduces interactions between the particles, in contrast to the BPF, where particles evolve independently from each other between resampling times.
Another interesting aspect of the FPF is the structure of the innovation term (Eq.~\eqref{Eqn:FPFparticles}, term in square brackets) that multiplies the gain function: 
the new observation $dY_t$ is compared to the arithmetic mean of the individual particle's estimate $h(Z_t^{(i)})$ and the population estimate $\bar{h}_t$. 
The function $\Omega$ comes from the conversion from Stratonovich form to It\^o form.

In practical implementations of the FPF, the main difficulty lies in the solution of the EL-BVP.
In multiple dimensions ($D>1$), na\"ive approaches to the EL-BVP turn out to be computationally expensive.
This problem has recently been addressed in \cite{Yang:2016ku} using a Galerkin approximation of the gain function.
We use the \emph{constant gain approximation} from that paper, which corresponds to a Galerkin approximation with coordinate functions as basis functions.
With this choice, the gain function is
\begin{equation}
K_{ij}(z,t)=\frac{1}{N}\sum_{k=1}^N \parenths{h_j(Z_t^{(k)})-\bar{h}_{t,j}}Z_{t,i}^{(k)}, \quad i,j=1,...,D,
\end{equation}
which is constant in $z$ (but not in time) and corresponds to the sample covariance matrix between the particle positions and the observation function $h$.
It can be shown that for a linear filtering problem, the constant gain approximation equals the Kalman gain in the limit of infinite number of particles \cite{Yang:2016ku}.
In addition, under this approximation the FPF gives rise to almost the same particle dynamics as in the ensemble Kalman-Bucy filter (EnKFBF) \cite{Bergemann:2012dx}, only differing slightly in the estimation of the covariance matrix from the observations.

\section{Results}
\label{results}

\subsection{Illustration of the curse of dimensionality for a linear problem}
We consider the linear $D$-dimensional problem
\begin{align}
\label{Eqn:state_lin}dX_t=-X_t dt+\sqrt{2}dW_t, \\
\label{Eqn:obs_lin}dY_t=2X_tdt+dV_t,
\end{align}
where the numerical factors are chosen in order to ensure a unit prior variance and time-constant and a (one-dimensional) mean squared error of the optimal filter $\mathds{E}[(\hat X_{i,t}-X_{i,t})^2]=1/2$ across all values of $D$.
For the linear system in Eqs.~(\ref{Eqn:state_lin},\ref{Eqn:obs_lin}), the bootstrap particle filter (BPF) is given by
\begin{align}
\label{Eqn:particlemotion_lin} dZ_t^{(i)} &= -Z_t^{(i)}dt+\sqrt{2}dB_t^{(i)}, \\
\label{Eqn:weightupdate_lin} m_t^{(i)}&=\frac{M_t^{(i)}}{\sum_{j=1}^NM_t^{(j)}}, \quad dM_t^{(i)}=2M_t^{(i)}Z_t^{(i)}\cdot dY_t,
\end{align}
where $i=1,...,N$ and $Z_0^{(i)}\sim\mathcal{N}(0,\mathds{1}_{D\times D})$ and $m_0^{(i)}=1/N$. 
We note that in principle all processes should carry an index that refers to the dimension $D$. 
However, in the interest of readability, we omit this index.

\subsubsection{The time-scale of weight decay shortens in higher dimensions} 
\begin{figure}[tbp]
\centering
\includegraphics[width=\textwidth]{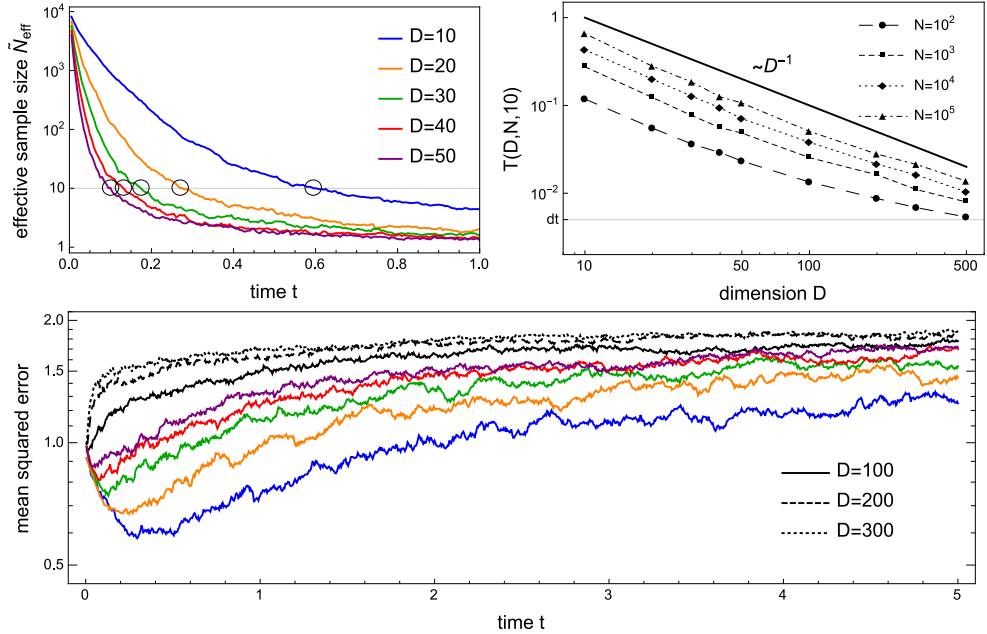}
\caption{\textbf{Top left:} The collapse of the effective sample size $\tilde N_{\text{eff},t}$ over time is shown for different dimensions $D$ of the state space and for an ensemble size $N=10^4$. The displayed time-course is an average of 100 independent trials. \textbf{Top right:} A plot of the stopping time $T(D,N,n)$ for $n=10$ as a function of dimension $D$ shows an approximate $D^{-1}$ scaling (thick black line). \textbf{Bottom:} The time evolution of the mean squared error shows a dip and a subsequent deterioration of performance due to weight decay for dimensions up to 50. The dip is shallower and the deterioration faster for higher dimensions. For even higher dimensions, the dip is no longer visible because the weight decay is too quick. All traces are averages of 100 independent trials with an ensemble size of $N=10^4$.}
\label{Fig:1}
\end{figure}

We study the effect of the number of dimensions $D$ on the time-scale on which the decay of $\tilde N_{\text{eff},t}$ occurs.
In Fig.~\ref{Fig:1} we show numerical results for the trial-averaged $\tilde N_{\text{eff},t}$, illustrating the decay of $\tilde N_{\text{eff},t}$ as a function of time for $N=10^4$ and $D=10,20,...,50$. 
The decay of $\tilde N_{\text{eff},t}$ critically limits the performance of the filter.
We measure the performance by the mean squared error (MSE) of the particle estimate of the hidden state
\begin{equation}
\text{MSE}_{t}=\frac{1}{D}\EE{\norm{X_t-\sum_{i=1}^Nm_t^{(i)}Z_t^{(i)}}^2},
\end{equation}
where the average is estimated numerically by averaging independent trials, as in the case of expected stopping time.
In our numerical experiment, the initial value is $\text{MSE}_0=1$ due to the initialization of the particles according to the prior distribution.
The MSE reaches an asymptotic value of two as the effective sample size goes to one.
As we show in Fig.~\ref{Fig:1} bottom, for low values of $D$ the MSE has a transient dip that is followed by a gradual increase towards the asymptotic value.
The dip becomes shallower and the increase becomes faster as the dimension $D$ increases.
For very high $D$, the dip disappears and the MSE increases immediately.

\begin{figure}[tbp]
\centering
\includegraphics[width=\textwidth]{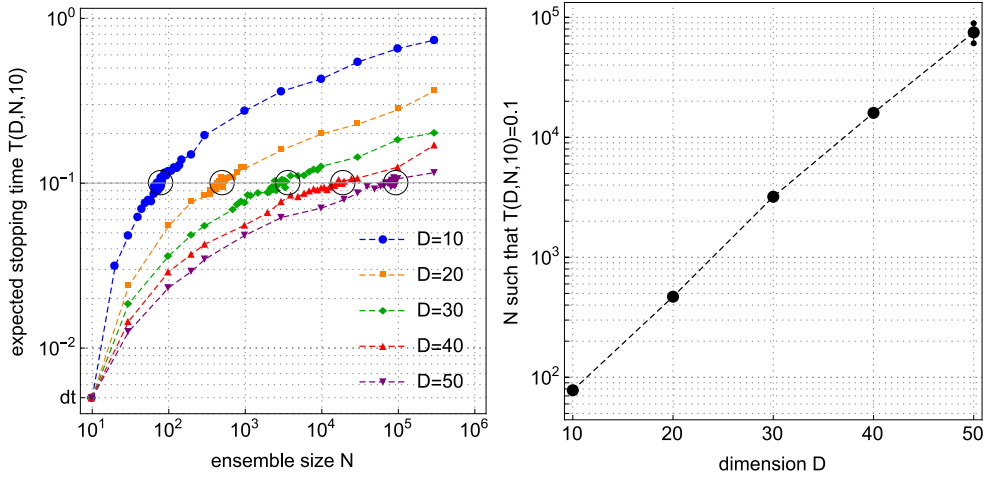}
\caption{\textbf{Left:} The time-scale $T$ of weight collapse (defined as the time in takes for $\tilde N_{\text{eff},t}$ to reach a value of 10, see black circles in Fig.~\ref{Fig:1} top left) as a function of ensemble size $N$ and dimension $D$. \textbf{Right:} The required ensemble size to achieve $T=0.1$ (see black circles in the left panel) increases exponentially with dimension $D$. }
\label{Fig:2}
\end{figure}

In order to quantify the time-scale of weight decay, we define the following expected stopping time, the mean first-passage time of $\tilde N_{\text{eff},t}$ through $n$,
\begin{equation}
T(D,N,n) = \EE{\inf\braces{t\geq 0 | \tilde N_{\text{eff},t}\leq n}}, \quad n\leq N,
\label{Eqn:stoppingtime}
\end{equation}
as a measure for the time-scale of weight degeneracy.
By definition, we have 
\begin{equation}
T(D,N,N)=0, \quad T(D,N,1)=\infty, \quad N\in\mathds{N}.
\end{equation}

We first give a rough analytical estimate of the dependence of $T(D,N,n)$ on $D$.
The time evolution of the normalized weight can be written as
\begin{equation}
dm^{(i)}_t=4m^{(i)}_t(Z^{(i)}_t-\dbar Z_t)\cdot(X_t-\dbar Z_t)dt+2m^{(i)}_t(Z^{(i)}_t-\dbar Z_t)\cdot dV_t,
\end{equation}
for the linear model.
From this we obtain
\begin{equation} 
\mathds{E}\brackets{d\parenths{m^{(i)}_t}^2}= \mathds{E}\brackets{4\parenths{m^{(i)}_t}^2\braces{2(Z^{(i)}_t-\dbar Z_t)\cdot(X_t-\dbar Z_t)+\norm{Z^{(i)}_t-\dbar Z_t}^2}dt},
\end{equation}
which consists of scalar products of $D$-dimensional vectors.
Each of the vector components is initially of order $1$, independently of $D$.
Even in the best of cases, i.e. if the particle positions are samples from the true posterior distribution, each of the vector components would be of the order of the true posterior standard deviation, which is equal to $1/\sqrt{2}$ independently of $D$.
Thus the initial magnitude of change of $\tilde N_{\text{eff},t}$, which is the inverse of the sum of the squared weights, is also proportional to $D$.
From this, a rough estimate of the scaling is $T(D,N,n)\propto D^{-1}$.
In Fig.~\ref{Fig:1} we numerically estimated $T$ by using trial averages, and we show the results as a function of $D$ for different values of $N$.
This confirms that the scaling is close to $D^{-1}$.

Next, we study the dependence of $T(D,N,n)$ on $N$, for which we rely on a numerical investigation.
The results are shown in Fig.~\ref{Fig:2} for $n=10$ and $D=10,20,...,50$.
We can see that as $D$ increases, $N$ has to increase exponentially in order to achieve a fixed $T$, or in other words $T(D,N,n)\propto \log N$.

\subsubsection{The effect of resampling}
\begin{figure}[tbp]
\centering
\includegraphics[width=\textwidth]{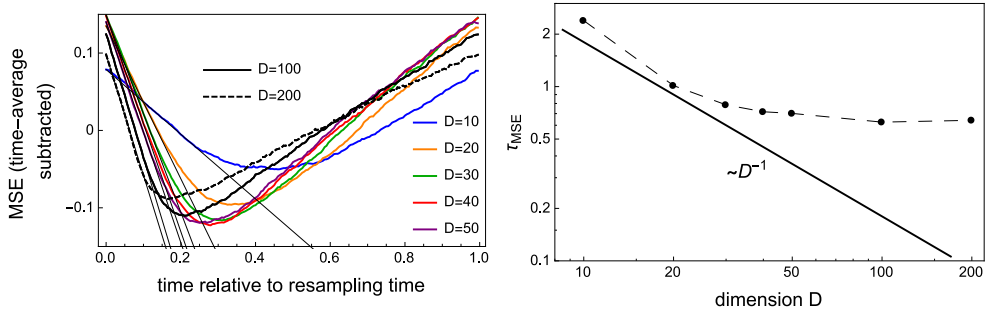}
\caption{The effect of resampling for an ensemble size of $N=10^4$ and a constant resampling interval of one time unit. 
\textbf{Left:} Resampling temporarily improves filter performance relative to the mean performance as shown by the dip in the time-course of MSE (with its time-average subtracted) after resampling. 
The thin straight lines show the slope of the initial decay. 
\textbf{Right:} We measure the slopes of the linear fits to the initial segment of the MSE time-course in the left panel. 
The characteristic time $\tau_{\text{MSE}}$, which is defined as the inverse absolute value of that slope, decays and stabilizes, whereas the time-scale of weight degeneracy continues to decay with $D^{-1}$ (thick black line, see also Fig.~\ref{Fig:1}).}
\label{Fig:3}
\end{figure}

If the criterion to resample the particles is $\tilde N_{\text{eff},t}=n_{\text{crit}}$, the rate at which resampling occurs is tied to the time-scale of weight degeneracy.
The immediate implication is that with a fixed ensemble size $N$, resampling occurs more frequently in higher dimension, with resampling rate roughly proportional to $D$.
After resampling, $\tilde N_{\text{eff},t}$ is reset to $N$, and since the resampled particles are located at positions where the likelihood of observations is high, the initial decay of $\tilde N_{\text{eff},t}$ is a bit slower than for an initialization from the prior.
However, since resampling does not add any new information about the true state, it cannot lead to an immediate performance increase.

The benefit of resampling is that particles with vanishing weights are discarded, and all computational efforts are expended for particles that are in an interesting region of state space.
It is therefore expected that resampling shows a delayed effect due to the diffusion of particles away from resampled positions.
For example, in the extreme case where $n_{\text{crit}}\approx 1$, all particles will typically be resampled at the same location.
The particles have to diffuse away from their initial position such that their empirical variance is of the same order of magnitude as the (true) posterior variance.
The time that is needed to reach such a state is related to the inverse of the squared diffusion coefficient and therefore independent of dimension.

In order to quantify the delay between the resampling time and the time where the MSE has maximally decreased due to resampling, we introduce the measure $\tau_{\text{MSE}}$, which is defined as the inverse absolute value of the slope of a linear fit to the initial portion of the MSE time-course (after resampling).
In the right panel of Fig.~\ref{Fig:3} we show $\tau_{\text{MSE}}$ as a function of dimension $D$.
The time-scale $\tau_{\text{MSE}}$ decreases with $D$, but it tends to a constant value for very high $D$.
This saturation is in sharp contrast to the time-scale of weight degeneracy, which continues to decay with $D^{-1}$.
The consequence is that in high dimensions, the weight decays during an interval which is much shorter than the delay required for the resampling to have a beneficial effect. 
Resampling is therefore ineffective in high dimensions, and it does not remove the need for exponentially large ensembles.

\subsection{Optimal proposals in continuous time}
In the literature on discrete time particle filters, optimal proposals for the particle motion have been shown to greatly reduce the required ensemble size.
However, proposals for discrete-time filters are not directly applicable to the continuous-time case.
For example, it is straightforward to show (c.f. Appendix~\ref{App:1}) that the optimal particle filter by \cite{Doucet:2000ui} collapses to a bootstrap particle filter as the time discretization step goes to zero.
Since the re-weighting of samples in the continuous-time particle filter depends on the mutual absolute continuity of the hidden state process and the particle process, the class of admissible particle motions is restricted.
Two diffusion processes are mutually absolutely continuous if and only if they differ by a pure drift term. 
The SDE for the particle motion can therefore differ from the hidden state SDE by at most a drift term $F_t$:
\begin{align}
\label{Eqn:particlemotion_alt} dZ_t^{(i)} = \brackets{f(Z_t^{(i)})+F^{(i)}_t}dt+g(Z_t^{(i)})dB_t^{(i)},
\end{align}
where $F_t^{(i)}$ must be a $\mathcal{F}_t^{Y}\vee\mathcal{F}_t^{Z^{(i)}}$-adapted $D$-dimensional process. The corresponding (unnormalized) weight will evolve as
\begin{align}
\label{Eqn:weightupdate_alt} dM_t^{(i)}=M_t^{(i)}\parenths{h(Z_t^{(i)})\cdot dY_t-F^{(i)}_t\cdot dZ_t^{(i)}}.
\end{align}
From this starting point, it could be interesting to formulate a stochastic control problem in order to choose the processes $F^{(i)}_t$ such as to minimize the effects of weight degeneracy.
There is some work in this direction \cite{VandenEijnden:2013iw} for continuous-discrete filters in the small observation noise regime, but we are not aware of any work on generalizing this to the continuously observed case (and for large observation noise).

It is also interesting to note that both the optimal proposal and the FPF have the aim of redirecting the particles based on the observations.
However, the particle motion of an importance sampling-based particle filter is more heavily constrained in order for the importance weights to exist.
Even the general form of Eq.~\eqref{Eqn:particlemotion_alt} is less general than the motion of the particles in the FPF, which include an explicit $dY_t$ term.
It remains an open problem to reconcile the two frameworks of particle filtering.

\subsection{Dimensionality-dependent scaling of the BPF vs. FPF}
\begin{figure}[tbp]
\centering
\includegraphics[width=\textwidth]{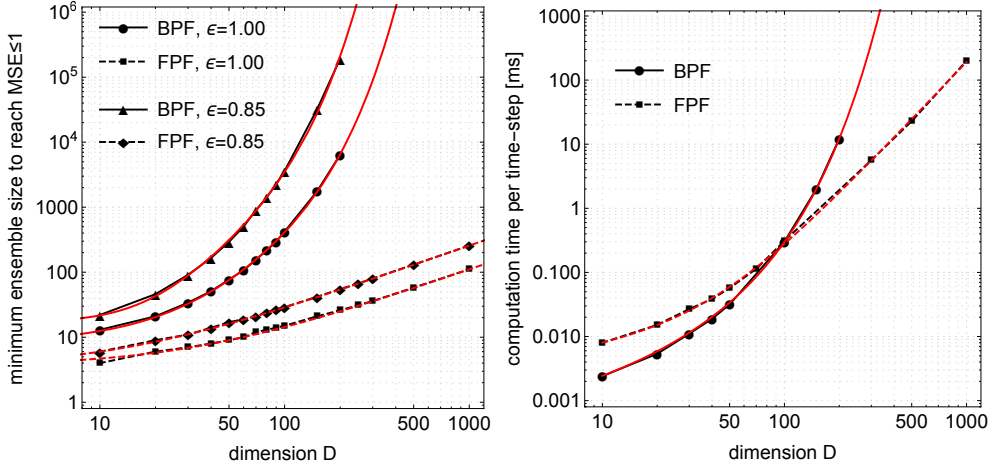}
\caption{\textbf{Left:} Comparison of the ensemble size $N_1(D)$ and $N_{0.85}(D)$ that is required in order to reach a performance of $\text{MSE}\leq 1$ and $\text{MSE}\leq 0.85$ respectively as a function of dimension $D$ for the Bootstrap Particle Filter (BPF) and the Feedback Particle Filter (FPF). The ensemble sizes are generally lower for the FPF and scale linearly with dimension, whereas the BPF requires an exponentially large ensemble. Fits are shown in red: $N^{\text{BPF}}_{1}(D)\approx  38.0\cdot e^{0.026 D} - 0.66\cdot D-30.2$, $N^{\text{BPF}}_{0.85}(D)\approx  68.1\cdot e^{0.040 D} - 2.93\cdot D-50.9$, $N^{\text{FPF}}_{1}(D) \approx 3.63+0.107 \cdot D$, and $N^{\text{FPF}}_{0.85}(D) \approx 3.46+0.253 \cdot D$. \textbf{Right:} Comparison of the run times of the two filtering algorithms for $\epsilon=1$. }
\label{Fig:4}
\end{figure}

We study the minimum ensemble size that is required in order to reach a certain performance (measured as mean squared error)
\begin{equation}
N_{\epsilon}(D) = \inf\braces{N\in\mathds{N} | \text{MSE}\leq \epsilon},
\label{Eqn:reqensemble}
\end{equation}
where MSE is the expected time-averaged mean-squared error defined as
\begin{equation}
\text{MSE}=\frac{1}{D t_1}\mathds{E}\int_0^{t_1}dt\begin{cases}\norm{X_t-\sum_{i=1}^Nm_t^{(i)}Z_t^{(i)}}^2 & \text{for the BPF} \\
\norm{X_t-\tfrac{1}{N}\sum_{i=1}^NZ_t^{(i)}}^2 & \text{for the FPF.}
\end{cases}
\end{equation}
In simulations, this integral is estimated using a Riemann sum over discrete time-steps, $t_1$ is chosen to be 5000 time units, for which we can drop the expectation because of the ergodicity of the process.
Particles of the BPF are resampled whenever $\tilde N_{\text{eff}}/N \leq 0.1$. 

For our particular linear toy model, the theoretical range of sensible values of $\epsilon$ is $0.5<\epsilon<2$, and we can set $N_{0.5}(D)=\infty$ and $N_{2}(D)=1$ because the performance of the optimal filter is set to produce an MSE of 0.5 and both filters go back to the prior for $N=1$, which yields an MSE of 2.
However, the practical range of $\epsilon$ is severely restricted by the runtime of the simulations, especially for the BPF.

The results for $\epsilon=1$ and $\epsilon=0.85$ are displayed in Fig.~\ref{Fig:4} left.
Already in 10 dimensions, the FPF starts out with a significant advantage, requiring only four particles vs. 13 for the BPF.
For larger dimensions, the ensemble size increases rapidly, reaching a value of 421 by D=100 and 6434 by D=200. 
Simulations for $D>200$ were too time-consuming to run.
Meanwhile, the FPF requires only 15 particles for $D=100$ and an estimated 25 particles for $D=200$. 
We ran simulations of the FPF up to $D=1000$, where it requires merely 111 particles.

A least-squares fit of the numerical data reveals an exponential scaling for the BPF and merely linear scaling for the FPF (see Fig.~\ref{Fig:4} caption).
The FPF thus requires roughly one additional particle with every increase of the dimension by ten.
In contrast, for the same increase in dimension, the BPF requires a factor of 1.3 more particles. 

In terms of computation time for a fixed performance of $\text{MSE}=1$, we find that the BPF requires a run time per time-step that scales exponentially with dimension (see Fig.~\ref{Fig:4} right).
In contrast, the FPF shows a cubic scaling, which is expected because the gain function scales quadratically with the ensemble size and the latter scales linearly with dimension.
Interestingly, the FPF requires more computation time for low dimensions, despite using fewer particles to achieve the same performance as the BPF. 
However, we did not heavily optimize our code for performance, and we expect that the runtime could be reduced to compete with the BPF in low dimensions using more careful programming.

\section{Discussion}
\label{discussion}
In this paper, we revisited the problem of curse of dimensionality (COD) in the standard particle filter.
We considered the case of the classical filtering problem with a continuous time index.
Even though the COD has been studied before, all the existing literature considers only one Bayesian update step and implies a discrete-time treatment.
Here, we closed this gap by studying the full dynamic nature of the problem in continuous time.

The discrete- and continuous-time particle filters have some important differences.
The class of possible proposal distributions is larger in discrete time, where it is only restricted in terms of practicality by virtue of tractability of the transition kernel.
In discrete time, it has been shown that even the optimal proposal distribution does not avoid the COD \cite{Snyder:2015dw}.
In continuous time, the law of the particle motion has to be absolutely continuous with respect to the law of the hidden state.
It is an open problem to show that there are nontrivial proposals that minimize the weight degeneracy.

There has been a general consensus that the problems of particle filters in high-dimensional problems result from importance sampling.
It has therefore been conjectured that a particle filter without importance weights could work efficiently in high dimensions.
Such a filter, the Feedback Particle Filter (FPF), has recently been proposed by \cite{Yang:hf}.
There have been other related approaches to filtering with unweighted particles, e.g. in \cite{Crisan:2010fm} and \cite{Crisan:2011tua}, Chapter 23, the Ensemble Kalman (Bucy) Filter \cite{Evensen:1994uq,Bergemann:2012dx} and the Neural Particle Filter in \cite{Kutschireiter:2017ip}.
Both the EnKBF and the Neural Particle Filter are mathematically similar to the FPF with constant gain approximation, and the latter differs only in terms of the structure of the innovation term.
Despite the promise of unweighted particle filters, their efficiency in high dimensions has not been thoroughly demonstrated so far.
This work complements the analyses in \cite{Yang:2016ku,Kutschireiter:2017ip} by providing numerical evidence that the FPF requires only a polynomial ensemble size and computation time as a function of dimensionality and thereby avoids the COD.
Future analytical work will hopefully shed light on the mechanisms that underly this finding.

Based on this important result, we want to draw the attention of researchers to the FPF and similar unweighted approaches.
Particle filters without importance weights are promising algorithms for solving very high-dimensional problems.
This opens up new perspectives in applied fields such as geophysics and meterology, especially numerical weather prediction, where data assimilation typically requires the solution of very high-dimensional filtering problems.

\appendix
\section{Continuum limit of the optimal proposal from \cite{Doucet:2000ui}}
\label{App:1}
In order to apply the optimal importance function that is given in \cite{Doucet:2000ui}, Example 3, we have to introduce a time-discretized version of our diffusion process:
\begin{align}
X_k &= X_{k-1}+f_{\text{cont}}(X_{k-1})dt+g(X_{k-1})\sqrt{dt}\epsilon_k,\\
Y_k-Y_{k-1} &= h(X_{k})dt+\sqrt{dt}\eta_k,
\end{align}
where $\epsilon_k$ and $\eta_k$ are multivariate Gaussian random variable with mean zero and unit covariance matrix.
For constant $g(x)=G$ and linear $h(x)=W x$, this corresponds to Eqs.~(9,10) in \cite{Doucet:2000ui} with $n_x=n_y=D$, $f(x)=x+f_{\text{cont}}(x) dt$, $\Sigma_v=dt GG^{\top}$, $C=W dt$, and $\Sigma_w=dt \mathds{1}_{D\times D}$, where the increment $Y_{k}-Y_{k-1}$ was identified with $y_k$ of \cite{Doucet:2000ui} in order to preserve the non-regressive form of $y$.
Taking Eqs.~(11,12) in \cite{Doucet:2000ui} and expressing the mean $m_k$ and covariance matrix $\Sigma$ of the conditional distribution $x_k|(x_{k-1},y_k)$ in terms of the quantities of the continuous-time model, we obtain
\begin{align}
m_k &= \Sigma\parenths{dt^{-1}G^{-\top}G^{-1}x_{k-1}+G^{-\top}G^{-1}f_{\text{cont}}(x_{k-1})+W^{\top}y_k},\\
\Sigma^{-1} &= dt^{-1}G^{-\top}G^{-1}+dt W^{\top}W.
\end{align}
For small $dt$, the covariance matrix can be approximated as $\Sigma\approx dt GG^{\top}$, and therefore we obtain
\begin{equation}
m_k = x_{k-1}+f_{\text{cont}}(x_{k-1})dt+dtGG^{\top}W^{\top}y_k.
\end{equation}
The last term, when expressed in terms of the continuous-time process, 
\begin{equation}
dZ^{(i)}_t=f_{\text{cont}}(Z^{(i)}_t)dt+G dW_t+GG^{\top}W^{\top}dY_tdt,
\end{equation}
constitutes a $dtdY_t$ term, which vanishes.
Therefore, the particle SDE reduces to
\begin{equation}
dZ^{(i)}_t=f_{\text{cont}}(Z^{(i)}_t)dt+G dW_t,
\end{equation}
which is equal to the prior.

\bibliographystyle{siamplain}
\bibliography{/Users/simone/Documents/uni/postdoc/library}
\end{document}